\documentclass[leqno,12pt]{amsart}
\makeatletter
\renewcommand\paragraph{\@startsection{paragraph}{4}{\z@}%
            {-2.5ex\@plus -1ex \@minus -.25ex}%
            {1.25ex \@plus .25ex}%
            {\normalfont\normalsize\bfseries}}
\makeatother
\makeatletter
\renewcommand{\@seccntformat}[1]{%
  \ifcsname format#1\endcsname
    \csname format#1\endcsname
  \else
    \csname the#1\endcsname\quad
  \fi
}
\makeatother

\DeclareFontFamily{U}{mathb}{}
\DeclareFontShape{U}{mathb}{m}{n}{
     <-5.5> mathb5
  <5.5-6.5> mathb6
  <6.5-7.5> mathb7
  <7.5-8.5> mathb8
  <8.5-9.5> mathb9
  <9.5-11>  mathb10
  <11->     mathb12
}{}
\DeclareSymbolFont{mathb}{U}{mathb}{m}{n}
\DeclareMathSymbol{\lefttorightarrow}{3}{mathb}{"FC}
\DeclareMathSymbol{\righttoleftarrow}{3}{mathb}{"FD}

\makeatletter
\DeclareRobustCommand{\looparrow}[1]{%
  \mathrel{\mathpalette\looparrow@{#1}}%
}
\newcommand{\looparrow@}[2]{\reflectbox{$\m@th#1#2$}}
\makeatother

\usepackage[margin=1in]{geometry}

\usepackage[
    colorlinks=true,
    linkcolor=blue,
    citecolor=blue,
    urlcolor=blue,
    bookmarks=true,         
    bookmarksnumbered=true, 
    bookmarkstype=toc      
]{hyperref}

\usepackage{bookmark}
\usepackage{amsmath}\usepackage{MnSymbol}\usepackage{phaistos}
\usepackage{amsfonts}\usepackage{verbatim}\usepackage{bbold}
\usepackage{graphicx}\usepackage{fontawesome}
\usepackage{marvosym}\usepackage{ifsym}
\usepackage{dictsym}\usepackage{bbm}
\usepackage{wasysym}\usepackage{dutchcal}
\usepackage{amstext}
\usepackage{amsbsy}
\usepackage{amsopn}
\usepackage{amsthm}\usepackage{color}
\usepackage{pigpen}
\newcommand{\CWrighttoleftarrow}{%
  \mathrel{\raisebox{\depth}{\scalebox{1}[-1]{$\righttoleftarrow$}}}%
}

\newtheorem{lmm}{Lemma}[section]

\newtheorem{prp}{Proposition}[section]

\newtheorem{thm}{Theorem}[section]

\theoremstyle{definition}

\newtheorem{expl}{Example}

\theoremstyle{remark}

\newtheorem*{thm*}{Theorem}
\def\proclaim#1{\vskip0.5em\noindent{\bf #1}\it}
\def\endproclaim{\vskip0.5em\par\noindent\rm}
\def\proclaim#1{\vskip0.5em\noindent{\bf #1}\it}
\def\endproclaim{\vskip0.5em\par\noindent\rm}
\def\demo#1{\vskip0.5em\noindent{\bf #1\ }}

\def\text#1{\mbox{#1}}
\def\flushpar{\par\noindent}

\def\tag#1{\eqno{(#1)}}

\def\e{\varepsilon}
\def\a{\alpha}
\def\b{\beta}

\def\d{\delta}

\def\D{\Delta}

\def\f{\flushpar}

\def\v{\varphi}

\def\om{\omega}

\def\B{\mathcal B}

\def\T{\widehat T}
\def\({\biggl(}
\def\){\biggr)}
\def\<{\langle}
\def\>{\rangle}
\def\bdy{\partial}\def\({\biggl(}
\def\){\biggr)}
\def\[{\biggl[}
\def\]{\biggr]}
\def\sbul{\f$\bullet\ \ \ $}
\def\pf{\smallskip\f{\tt Proof}\ \ \ \ }\def\sms{\smallskip\f}
\def\with respect to {\text{with respect to}} 
\def\im{\text{\tt\small  Im}\,}
\def\ttau{\widehat{\tau}}

\let\oldtocsection=\tocsection

\let\oldtocsubsection=\tocsubsection

\let\oldtocsubsubsection=\tocsubsubsection

\renewcommand{\tocsection}[2]{\hspace{0em}\oldtocsection{#1}{#2}}
\renewcommand{\tocsubsection}[2]{\hspace{1em}\oldtocsubsection{#1}{#2}}
\renewcommand{\tocsubsubsection}[2]{\hspace{2em}\oldtocsubsubsection{#1}{#2}}
\begin{document}
\title{Dual mixing and quasi-mixing  of inner function restrictions   }

\author{Jon Aaronson}
\address[]{School of Math. Sciences, Tel Aviv University,
69978 Tel Aviv, Israel.}
\email{aaro@tau.ac.il}
\begin{abstract}Dual quasi-mixing and dual mixing of a measure preserving transformation are
uniform individual ratio limit properties of the transfer operator which entail the quasi-mixing and mixing properties of Krickeberg (respectively). The real restrictions of  normalised, parabolic, non-M\"obius inner functions of the upper half plane preserve Lebesgue measure and are dual quasi mixing; such a restriction being dual mixing iff it is exact 
and singular dual quasi-mixing otherwise.
 \end{abstract}
\thanks{\copyright 2026.}
\maketitle\markboth{\copyright Jon. Aaronson}{mixing of inner functions}
\tableofcontents\setcounter{tocdepth}{3}
\section{ Introduction}
\
\subsection{Overview}
\

Krickeberg (\cite{Krickeberg,Krickeberg-recent})  defined  notions of topological  mixing and quasi-mixing for infinite measure preserving transformations (after   an example of E.Hopf (\cite[\S17]{Hop1}))\ \footnote{called``Hopf-Krickeberg mixing/quasi-mixing'' respectively  in the sequel}.
\

A stronger, transfer operator  version of Hopf-Krickeberg mixing (which we call ``dual mixing'') was established for certain ``intermittent'' interval maps in \cite{Thaler2000}.
\

This paper considers ``dual mixing'' and ``dual quasi-mixing'' of ``normalised, parabolic inner function restrictions'' which  provide new  examples of Hopf-Krickeberg quasi-mixing transformations.

\

We conclude this overview with two illustrative examples, both real restrictions of rational, normalised,  parabolic,  inner functions of the upper half plane (as in \S\ref{Inner}).
\

Let $(\Bbb R,m)$ be the real line equipped with Lebesgue measure.
\expl\label{Boole}

Let $T:\Bbb R\CWrighttoleftarrow$ be {\it Boole's transformation}
as in \cite{Ad-W}:
$$T(x):=x-\tfrac1x$$ and set $u_n=\tfrac1{\pi\sqrt{2n}}$,
then $(\Bbb R,m,T)$ is a measure preserving transformation and
\begin{align*}\tag*{\Coffeecup}\label{Coffeecup}\tfrac1{u_n}\T^nf\xrightarrow[n\to\infty]{\text{\scriptsize u.c.s.}}\int_\Bbb R fdm\ \forall\ f\in C_C(\Bbb R)
\end{align*}
where $C_C(\Bbb R):=\{\text{\tt\small continuous functions on $\Bbb R$  with compact support}\}$.
\

Here $\T:L^1(m)\CWrighttoleftarrow$
is the {\it transfer operator} as defined  by\ \  \ref{faPaw}\ \ on p.\pageref{faPaw}  and u.c.s means ``uniformly on compact subsets''.
\

In other words, $(\Bbb R,m,T)$ is  {\tt dual mixing} as defined by \ref{dsmilitary}
(on p.\pageref{dsmilitary}).
\

\expl\label{Boole+1}
\

Define $\tau:\Bbb R\CWrighttoleftarrow$ by
$$\tau(x):=T(x)+1=x-\tfrac1x+1,$$
then $(\Bbb R,m,\tau)$ is a measure preserving transformation and
there is a singular, Borel measure $\rho\nequiv 0$ on $\Bbb R$ which is  locally finite (in the sense that compact sets have finite $\rho$-measure)   so that with  $v_n=\tfrac1{n^2}$,
\begin{align*}\tag*{\Football}\label{Football}\tfrac1{v_n}\ttau^nf\xrightarrow[n\to\infty]{\text{\scriptsize u.c.s.}}\int_\Bbb R fd\rho\ \forall\ f\in C_C(\Bbb R).
\end{align*}
\

In other words, $(\Bbb R,m,T)$ is  {\tt singular dual quasi-mixing} as  defined by \
 \ \ref{dsliterary}\ \ (below).\ \

Proposition\ \  \ref{PS}\ \  on p.\pageref{PS}   establishes\ \  \ref{Coffeecup}\ \  and \ \ \ref{Football} .
\subsection{Dual  mixing and quasi-mixing}\label{DQM}
\

Let $X$ be a locally compact, Polish  space, let $m$ be a locally finite, Borel measure on $X$
and let $(X,m,T)$ be a measure preserving transformation
with {\it transfer operator} $\T:L^1(m)\CWrighttoleftarrow$ defined by
\begin{align*}\tag*{\faPaw}\label{faPaw}
\int_X\T g\cdot hdm=\int_Xg\cdot h\circ Tdm\ \text{for}\ g\in L^1(m),\ h\in
L^\infty(m).
\end{align*}

\

We'll call
 $(X,m,T)$ {\it dual   mixing}
 if  $\exists$  $u_n>0$ so that
 \begin{align*}\tag*{\dsmilitary}\label{dsmilitary}
  \tfrac1{u_n}&\T^n f(x)\xrightarrow[n\to\infty]{\text{\scriptsize u.c.s}}\ m(f):=\int_Xfdm\  \forall\ f\in C_C(\Bbb R).
 \end{align*}
See \cite[Chapter 3]{Orey-book}, and \cite{Thaler2000}.

.

More generally,
for
 $\mu\nequiv 0$  a locally finite, Borel measure  on $X$,
we'll call
 $(X,m,T)$ {\it $\mu$-dual   quasi-mixing}
 if  $\exists$  $u_n>0$ so that
 \begin{align*}\tag*{\dsliterary}\label{dsliterary}
  \tfrac1{u_n}&\T^n f(x)\xrightarrow[n\to\infty]{\text{\scriptsize u.c.s}}\mu(f):=\int_Xfd\mu\ \  \forall\ f\in C_C(\Bbb R).
 \end{align*}
In this case,  $(X,m,T)$ is {\it dual mixing}  if $\mu=m$;
and we'll call  $(X,m,T)$   {\it singular dual quasi-mixing} if $\mu\perp m$.

\

\f{\bf Remark}
\

\ It is standard to extend\ \  \ref{dsliterary}\ \  to the class of bounded,
$\mu$-integrable functions with compact support which are continuous at $\mu$-a.e. point of $\Bbb R$, whence  if $(X,m,T)$ is   $\mu$-dual quasi-mixing,   then
\begin{align*}m(A\cap T^{-n}B)&=\int_B\T^n 1_Adm\sim u_nm(B)\mu(A)\ \text{as}\ n\to\infty\\ & \forall\ A\ \text{bounded},\ \mu(\bdy A)=0\ \&\ B\in\B(X),\ B\ \text{bounded}.
\end{align*}
In particular, $(X,m,T)$ is $(\mu,m)$-Krickeberg quasi-mixing as in\, \ref{dsbiological}\, below.
\

Dual mixing for certain intermittent interval maps (including Boole's transformation as in\ \ \ref{Coffeecup}) was established in \cite{Thaler2000}. For more examples, see\ \cite{Melbourne-Terhesiu2012}.
\subsection{Hopf-Krickeberg  mixing and  quasi-mixing}\
\

Let $X$ be a locally compact, Polish  space and let $m$ be a locally finite, Borel measure on $X$.

A measure preserving transformation $(X,m,T)$ is called
\

 \f {\it Hopf-Krickeberg $(\mu_1,\mu_2)$-quasi-mixing} (abbr. HK-$(\mu_1,\mu_2)$-quasi-mixing)
  where $\mu_1,\mu_2\nequiv 0$ are locally finite, Borel measures on $X$;  if $\exists\ u_n>0$ so that
\begin{align*}\tag*{\dsbiological}\label{dsbiological}
\tfrac{m(A\cap T^{-n}B)}{u_n}\xrightarrow[n\to\infty]{}\mu_1(A)\mu_2(B)
\ \forall\ A\in\mathcal{R}_{\mu_1},\ B\in\mathcal{R}_{\mu_2}
 \end{align*}
 where $\mathcal{R}_\mu:=\{A\subset X:\ \text{precompact},\ \mu(\bdy A)=0\}$;
  and $(X,m,T)$ is called
 {\it Hopf-Krickeberg  mixing} (abbr. HK-mixing) if $(X,m,T)$ is\ \  HK-$(m,m)$-{quasi-mixing}.
\

For Markov shifts, HK-mixing is characterised by Orey's {\tt strong ratio limit property} (abbr. {\tt SRLP})  as in \cite{Orey}. See \cite{Krickeberg,Krickeberg-recent,Papangelou}.
\

In addition, a variety of smooth dynamical systems are now known to be HK-mixing. See
e.g. \cite{Lenci,DN,FR,PT}.
\

A characterisation of   HK-quasi-mixing Markov shifts with discrete state spaces in terms of
Pruitt's generalized {\tt SRLP} (as in \cite{Pruitt-SRLP})
is given by \cite[Theorem 2.4]{Papangelou}. Examples of
HK-quasi-mixing Markov shifts with continuous state space are also presented in \cite{Papangelou}.
  \
\

\subsection{Ergodic properties}
\

\subsubsection*{Dual quasi-mixing transformations }
\

A  dual   quasi-mixing measure preserving transformation $(X,m,T)$ of a probability space is {\it exact} in the sense that $\bigcap_{n\ge 0} T^{-n}\B(X)\overset{m}=\{\emptyset,X\}$.
\

By Hopf's decomposition results (see e.g. \cite[Proposition 1.3.1]{A1}), if $(X,m,T)$ is     $\mu$-dual   quasi-mixing, then $(X,m,T)$ is conservative iff $\sum_{n\ge 1}u_n=\infty$.
\

In this case, by the ratio ergodic theorem, almost every {\tt ergodic component}
$(X,m_\om,T)$ (of the { ergodic decomposition} of $(X,m,T)$ as in e.g. \cite[2.2.9]{A1}) is
{\tt\small pointwise dual ergodic} (as in \cite[\S3.7]{A1}) with $m_\om=\mu$, whence
$\mu=m$ and $(X,m,T)$ is  dual mixing which in turn entails {\tt\small conditional rational weak mixing}. See  \cite{A-T}.
\

\subsubsection*{HK-quasi-mixing transformations}
\

As shown in \cite{Krickeberg,Krickeberg-recent,Papangelou}, a  probability preserving transformation is HK-quasi-mixing iff it is strongly mixing.
\

A  {\tt rationally ergodic}  measure preserving transformation $(X,m,T)$ (see e.g. \cite[\S3.3]{A1}) which is $(\mu_1,\mu_2)$-HK-quasi-mixing as in\ \  \ref{dsbiological}\ \   on p.\ \pageref{dsbiological}\ \ with $\mathcal{R}_{\mu_i}\cap R(T)\ne\emptyset$ ($i=1,2$) is necessarily HK-mixing and hence by
\cite[Lemma C]{RWM} {\tt rationally weakly mixing}.

However, we do not know whether a general, conservative, HK-mixing transformation preserving an infinite measure is necessarily ergodic.

\subsection{Inner functions of the upper half plane}\label{Inner}

\

An analytic endomorphism $T:\Bbb R^{2+}\CWrighttoleftarrow$ is called an {\it inner function} of the {\it upper half plane}
$\Bbb R^{2+}:=\{x+iy:\ x,y\in\Bbb R,\ y>0\}$ if
$$\text{for $m$-a.e.}\ x\in\Bbb R,\ \exists\ \lim_{y\to 0+}T(x+iy)=:T(x)\in\Bbb R.$$

An inner function $T:\Bbb R^{2+}\CWrighttoleftarrow$ is necessarily of form

 \begin{align*}\tag*{\dsrailways}\label{dsrailways}T(z)=\a_T z+\b_T+\int_\Bbb R\tfrac{1+tz}{t-z}d\mu_T(t)\ \ (z\in\Bbb R^{2+}).
   \end{align*}
where $\a_T\ge 0,\ \b_T\in\Bbb R$ and $\mu_T$ is a  singular measure  on $\Bbb R$
(i.e. $\mu_T\perp m$).

 \

 The measurable mapping $T:\Bbb R\CWrighttoleftarrow$ is aka the {\it real restriction} of the inner function.

 For $\om\in\Bbb R^{2+},$ define
 $\v_\om(t):=\text{\tt\small Im}\,\frac1{\pi(t-\om)}\ \&\ P_\om\in\mathcal{P}(\Bbb R),\
 dP_\om:=\v_\om dm$.

\thm {\rm Nordgren's Theorem\ \cite{Nor} }\label{Nor}
\

For $T:\Bbb R^{2+}\CWrighttoleftarrow$ inner, as in \ref{dsrailways},  $(\Bbb R,m,T)$ is a nonsingular transformation with
$\T\v_\om=\v_{T\om}$; equivalently
\begin{align*}\tag*{\Bat}\label{Bat}P_\om\circ T^{-1}=P_{T(\om)}\ \forall\ \om\in\Bbb R^{2+}.\end{align*}\endproclaim

\

\f{\bf Remarks}

\sbul The case of Theorem \ref{Nor}  for rational inner functions was established
 in \cite[p. 787]{Boo} (see also \cite[\S8]{Glaisher}).

 \sbul For a converse to Nordgren's theorem,  see  \cite[Theorem]{Let}.

 \sbul As in \cite[Corollary]{Let},  if $\a_T>0$, then
\begin{align*}\tag*{\faAnchor}\label{faAnchor} m\circ T^{-1}=  \tfrac1{\a_T}m.\end{align*} \endproclaim

For $T:\Bbb R^{2+}\CWrighttoleftarrow$, we'll call the nonsingular transformation $(\Bbb R,m,T)$ the {\it restriction} of $T$.

\

\thm {\rm Denjoy-Wolff Theorem}\label{DW}
\

Let $T:\Bbb R^{2+}\CWrighttoleftarrow$ be a non-M\"obius inner function, then $\exists\ \mathfrak{d}=\mathfrak{d}_T\in\Bbb R^{2+}\cup\Bbb R\cup\{\infty\}$ so that
$$T^nz\xrightarrow[n\to\infty]{}\mathfrak{d}\
\forall\ z\in\Bbb R^{2+}.$$
\sbul If $\mathfrak{d}\in\Bbb R^{2+}$, then
$T(\mathfrak{d})=\mathfrak{d}$ and $|T'(\mathfrak{d})|<1$.
\sbul$\mathfrak{d}=\infty$ iff $\a_T\ge 1$.
\sbul  If $\mathfrak{d}\in\Bbb R$, then $\a_{T_{\mathfrak{d}}}\ge 1$ where $T_{\mathfrak{d}}:=\phi_{\mathfrak{d}_T}\circ T\circ\phi_{-\mathfrak{d}_T}$  with $\phi_t(z):=\tfrac{1+tz}{t-z}$.
\endproclaim

\

\f{\bf Denjoy-Wolff points}

The point $\mathfrak{d}_T$ is called the {\it Denjoy-Wolff point of $T$}.
\

The inner function $T$ is called
 {\it elliptic} if $\mathfrak{d}\in\Bbb R^{2+}$;
 {\it parabolic}  if $\mathfrak{d}\in\Bbb R\cup\{\infty\}$ and $\a_{T_{\mathfrak{d}}}=1$ (with $T_\infty:=T$); and
 {\it hyperbolic} if $\mathfrak{d}\in\Bbb R\cup\{\infty\}$ and $\a_{T_{\mathfrak{d}}}>1$ .

 We'll also call a non-elliptic inner function $T$ {\it normalised} if $\mathfrak{d}_T=\infty$. A non-elliptic inner function $T$ is M\"obius-conjugate
 to the  normalised inner function $T_{\mathfrak{d}}$ as above.
 \

 If $T:\Bbb R^{2+}\CWrighttoleftarrow$ is parabolic and normalised, then by \ref{faAnchor}, $(X,m,T)$ is a measure preserving transformation.

\thm {\rm Pommerenke's limit theorem\ (\cite[Theorem 1]{Pom})}\label{Pom}
\

Suppose that
$T:\Bbb R^{2+}\CWrighttoleftarrow$ is normalised, inner
and not M\"obius. Let $T^n(i)=a_n+ib_n$, then
\f{\rm (i)}\ $\exists\ \Pi:\Bbb R^{2+}\CWrighttoleftarrow$ analytic, so that
\begin{align*}
\tag*{\dsheraldical}\label{dsheraldical}&\tfrac{T^n-a_n}{b_n}\xrightarrow[n\to\infty]{\text{\tiny u.c.s.}}\Pi \ \text{with}\  \Pi\circ T\equiv a\Pi+b\ \ \text{where}\  b\in\Bbb R, a\ge\a_T.
\end{align*}
\

\f{\rm (ii)}\  If $\a_T=1$, then $a=1$; and if $T\ne\text{\tt Id}$, then
\begin{align*}\tag*{\dsagricultural}\label{dsagricultural}
\Pi\equiv i\ \ &\Leftrightarrow\ \Pi\circ T\equiv \Pi\ \Leftrightarrow\ \<T^n(z),T^n(\om)\>\xrightarrow[n\to\infty]{} 0
\ \forall\
z,\om\in\Bbb R^{2+}\\ &\text{where}\ \ \<z,\om\>:=|\tfrac{z-\om}{z-\overline{\om}}|.
\end{align*}
\f{\rm (iii)}\ \cite[Remark 1]{Pom}\ \ If $\a_T=1\ \&\ \Pi\nequiv i$, then $T^n(i)=a_n+ib_n\to\infty$ {\tt tangentially} in the sense that
$\D_n:=(\tfrac{a_n}{b_n})^2\to\infty$.
\endproclaim
We'll call the function $\Pi=\Pi_T:\Bbb R^{2+}\CWrighttoleftarrow$ the {\it Pommerenke limit of $T$} and the inner function $T$ {\it exact} if the Pommerenke limit $\Pi\equiv i$. This corresponds to the exacness of the restriction $(\Bbb R,m,T)$ (see \S\ref{erg}).

Hamilton showed (\cite[Theorem P]{Hamilton}) that the Pommerenke limit is either constant or  inner.

Since $V:=\text{\tt Im}\,\Pi$ is positive,  harmonic $\&\ V(i)=1,\ \exists\
\varpi_\Pi\in\mathcal{P}(\widehat{\Bbb R})$ so that
\begin{align*}\tag*{\faLock}\label{faLock}V(z)=\int_{\widehat{\Bbb R}}\tfrac{\v_z}{\v_i}d\varpi_\Pi=\varpi_\Pi(\{\infty\})\text{\tt\small Im}z+\int_\Bbb R\tfrac{\v_z}{\v_i}d\varpi_\Pi.
 \end{align*}
\

We'll call the measure $\varpi_\Pi$ the {\it representing measure} of $\Pi$.
Thus (\cite[Theorem P]{Hamilton})
either $\varpi_\Pi=P_i$ (when $V\equiv i$) or $\varpi_\Pi\perp P_i$ (otherwise).
\

Note that $\text{\tt spt}\varpi_\Pi=\{\infty\}$ $\iff$ $\text{\tt\small Im}\,\,\Pi(z)=\text{\tt\small Im}\,z$   which is  the case if and only if $T$ is
M\"obius. Thus,\  $T$ non-M\"obius entails $\varpi_{\Pi}(\Bbb R)>0$.
\subsection{Ergodic theory of inner function restrictions}\label{erg}
\

It follows from \ref{Bat}\ \ (p.\pageref{Bat})\     that
\begin{align*}\tag*{\dstechnical}\label{dstechnical}\int_{\Bbb R}|\v_z-\v_\om|dm={4\over\pi}\sin^{-1}\<z,\om\>.
\end{align*} See e.g. \cite[Lemma 6.1.4]{A1}.

\

By \cite{Lin-mixing71}, (see also e.g. \cite[Theorem 1.3.3]{A1}), the nonsingular transformation $(\Bbb R,m,T)$ is {\it exact}
in the sense that $\bigcap_{n\ge 0}T^{-n}\B(X)\overset{m}=\{\emptyset,\Bbb R\}$ iff
\begin{align*}\tag*{\dsmathematical}\label{dsmathematical}\|\T^nu\|_{L^1(m)}\xrightarrow[n\to\infty]{}0\ \forall\ u\in L^1(m),\ \int_Xudm=0;
 \end{align*}
whereas by approximation, using\ \  \ref{Bat}\ \ (p.\pageref{Bat});\ \ $\&$\  \
\ref{dstechnical}, we have that the  inner function restriction $(\Bbb R,m,T)$ saisfies
\ref{dsmathematical} iff
\begin{align*}&
\<T^nz,T^n\om\>\xrightarrow[n\to\infty]{}0\ \forall\ z,\om\in\Bbb R^{2+}.
\end{align*}
\

By Pommerenke's theorem \ref{Pom} (p.\pageref{Pom}), in case\ $(\Bbb R,m,T)$ is not exact, the Pommerenke limit is not constant and there is a non-constant $T$-invariant,  bounded analytic function on $\Bbb R^{2+}$, whence $(\Bbb R,m,T)$ is not ergodic.

\sms By Hopf's criterion (see e.g. \cite[Proposition 1.3.1]{A1}), $(\Bbb R,m,T)$ conservative iff $\sum_{n\ge 1}\T^n\v_i=\infty$ a.e.;  equivalently, $\sum_{n\ge 1}\tfrac1{b_n}=\infty$ with $b_n:=\tfrac{|T^n(i)|^2}{\text{\tt\small Im}\,T^n(i)}$.
\

As in \cite{Neu}  (see also \cite{A-inner78}\ $\&$\ \cite[Proposition 6.1.7]{A1}), if a non-M\"obius inner function is conservative, it is ergodic, hence exact as above.
\

However, as in \cite{A-inner81} (also \cite[Theorem 6.4.8]{A1}) there is a dissipative exact inner function restriction.

\section{Dual   mixing and   quasi-mixing of inner function restrictions}

 This section is devoted to  proving the

\proclaim{Main Theorem}\hypertarget{target:Main}{}
\

Let $T:\Bbb R^{2+}\CWrighttoleftarrow$ be a normalised, parabolic, non-M\"obius, inner function.
\sms  If $T$ is exact, then  $(\Bbb R,m,T)$ is   dual mixing.
\

\sms If $T$ is not exact, then $(\Bbb R,m,T)$ is
 singular dual quasi-mixing.
\endproclaim

\subsection{Uniform weak convergence of measures}
 \

Let $\mathcal{M}(Z)$ denote the space of finite Borel measures on the Polish space $Z$ and let
$\mathcal{P}(Z)$ denote the Borel probabilities on $Z$.
\

By Riesz's representation theorem, for $Z$ compact, $\mathcal{M}(Z)=C(Z)_+^*$ and by Helly's theorem $\mathcal{P}(Z)$ is weak * compact.
\

Fix   compact metric spaces $Y\ \&\ Z$. We'll call a function  $\mu:Y\to\mathcal{M}(Z)$  {\it weakly continuous} if
 $x\mapsto\int_Xfd\mu_x$ is continuous ($Y\to\Bbb R$) for each $f\in C(Z)$.
 Denote
 $$C(Y,\mathcal{M}(Z)):=\{\mu:Y\to\mathcal{M}(Z):\ \mu\ \text{weakly continuous}\}$$ considered with the topology of convergence
 $$\mu_{n}\xrightarrow[n\to\infty]{C(Y,\mathcal{M}(Z))}\nu\ \iff\ \sup_{x\in Y}\,|\mu_{n,x}(f)-\nu_x(f)|\xrightarrow[n\to\infty]{}0\ \forall\ f\in C(Z).$$

For $X$ be a locally compact, non-compact, Polish space,  let
  $$C_{\text{\tt\tiny lim}}(X):=\{f\in C(X):\ \exists\ \lim_{x\to\infty}f(x)\in\Bbb R\},$$
  then  $C_{\text{\tt\tiny lim}}(X)\ =\{f|_X:\ f\in C(\widehat{X})\}$ where
   $\widehat{X}:=X\cup\{\infty\}$ is the {\it one-point compactification} of $X$.

 \prp\label{precompacthat}
\

 Let $(Y.d)$ be a compact metric space,  then
 \

 $\mathcal{K}\subset C(Y,\mathcal{M}(\widehat{X}))$ is  precompact iff for some $J\subset C_{\text{\tt\tiny lim}}(X)_+$ with dense linear span,
 \begin{align*}\tag*{\faLeaf}\label{faLeaf}
 \om_f(h)&:=\sup\,\{|\mu_x(f)-\mu_y(f)|:\ x, y\in Y,\
 d(x,y)\le h,\ \mu\in \mathcal{K}\}\\ &\xrightarrow[h\to 0+]{}0\ \forall\ f\in J.
 \end{align*}
\endproclaim
\demo{Proof} \ Precompactness implies
\, \ref{faLeaf}\, by the Arzela-Ascoli theorem.
\

Conversely, suppose that $\mathcal{K}$ satisfies \ref{faLeaf}.
\

We claim first that $\mathcal{K}$ is  {\it uniformly bounded} in the sense  that
\begin{align*}\tag*{\faTree}\label{faTree}C:=\sup\,\{\mu_x(\mathbbm{1}):\ \mu\in\mathcal{K}, \  x\in Y\}<\infty\ \text{where}\ \mathbbm{1}(x)=1\ \forall\ x\in X. \end{align*}
If not,   $\exists\ n_k\uparrow\infty,\ x_k\in Y$ so that
$\mu_{n_k,x_k}(\mathbbm{1})\ge k$.
\

Fix $f\in \text{\tt\small Span}\,J$ so that $f(x)\ge c>0\ \forall\ x\in\widehat{X}$.
\

By \ref{faLeaf} and the Arzela-Ascoli theorem,
$\exists\ k_\ell\to\infty\ \&\ g\in C(Y)$ so that
$$\varlimsup_{\ell\to\infty}\mu_{n_{k_\ell},x_{k_\ell}}(f)\le C:=\sup_{y\in Y}|g(y)|$$
entailing for large $\ell$,
$$3C<k_\ell \le \mu_{n_{k_\ell} ,x_{k_\ell} }(\mathbbm{1})\le 2\mu_{n_{k_\ell} ,x_{k_\ell} }(f)\lesssim 2C.
\ \XBox\ \CheckedBox\ \text{\ref{faTree}}$$
To see precompactness of $\mathcal{K}$ fix an infinite subset $\mathcal{K}_0\subset\mathcal{K}$.
\

By \ref{faLeaf},
$\exists\ \mu_k\in\mathcal{K}_0$ so that
$$\mu_{k,x}(f)\xrightarrow[k\to\infty]{\text{\tiny uniformly}}G(f)(x)\ \forall\ f\in \text{\tt\small Span}\,J$$
  where $f\mapsto G(f)$ is  positive, linear
$\text{\tt\small Span}\,J\to C(Y)$.
\

Evidently $\|G(f)\|_{C(Y)}\le C\|f\|_{C_{\text{\tt\tiny lim}}(X)}$
with $C$ as in \ref{faTree}
and so $G$ extends by uniform
 approximation to a  positive, linear map $G:C_{\text{\tt\tiny lim}}(X)\to C(Y)$
satisfying $\|G(f)\|_{C(Y)}\le C\|f\|_{C_{\text{\tt\tiny lim}}(X)}$.
\

By the Riesz representation theorem $\exists\ \nu:Y\to\mathcal{M}(\widehat{X})$ so that
$$G(f)(x)=\int_Xfd\nu_x\ \forall\ f\in C_{\text{\tt\tiny lim}}(X),$$
\

Next,  fix  $f\in C_{\text{\tt\tiny lim}}(X)\ \&\ f_n\in\text{\tt\small Span}\,J$
with $\|f_n-f\|_{C_{\text{\tt\tiny lim}}(X)}=:\e_n
\xrightarrow[n\to\infty]{}0$. This entails $\|G(f_n)-G(f)\|_{C(Y)}\le C\e_n$, whence,
 $\forall\ k,n\ge 1$,
\begin{align*}|\mu_{k,x}&(f)-\int_{\widehat{X}}fd\nu_x| =|\mu_{k,x}(f)-G(f)(x)|\\ &\le
|\mu_{k,x}(f_n)-G(f_n)(x)|+|\mu_{k,x}(f)-\mu_{k,x}(f_n)|+|G(f)(x)-
G(f_n)(x)|\\ &\le |\mu_{k,x}(f_n)-G(f_n)|+2C\e_n.
\end{align*}
Given $\mathcal{E}>0$, fix $N=N_\mathcal{E}\ge 1$ so that
$\e_N<\tfrac{\mathcal{E}}{4C}$ and then choose $K=K_\mathcal{E}\ge 1$ so that
$$|\mu_{k,x}(f_N)-G(f_N)|<\tfrac{\mathcal{E}}{2C}\ \forall\ x\in Y,\ k\ge K$$
entailing
$$|\mu_{k,x}(f)-G(f)(x)|\le\mathcal{E}\ \forall\  x\in Y,\ k\ge K.$$
Thus $\nu\in C(Y,\mathcal{M}(\widehat{X}))$ and  $\mu_{n}\xrightarrow[n\to\infty]{C(Y,\mathcal{M}(\widehat{X}))}\nu$.\ \CheckedBox

\

\

For $Y\subset\Bbb R$ a bounded interval, define the {\it Lipschitz norm} of $F:Y\to\Bbb R$ by
$$\|F\|_{\text{\tt\tiny Lip(Y)}}:=\|F\|_{L^\infty(Y,m)}+\text{\tt Lip}_Y(F)$$ where
$$\text{\tt Lip}_Y(F):=\sup_{x,y\in Y,\ x\ne y}\tfrac{|F(x)-F(y)|}{|x-y|}.$$
Note that if $\text{\tt Lip}_Y(F)<\infty$, then $F$ is absolutely continuous on $Y$ and $\text{\tt Lip}_Y(F)=\|F'\|_{L^\infty(Y,m)}.$
\lmm\label{transferlem}
\

  Let $T:\Bbb R^{2+}\CWrighttoleftarrow$ be a normalised, parabolic, non-M\"obius, inner function with Pommerenke limit $\Pi=U+iV$, then for $z\in\Bbb R^{2+}$ and $Y\subset\Bbb R$
  a bounded interval,

\begin{align*}&\tag*{{\large\Pointinghand}}\label{Pointinghand}
F_{n,z}:=\tfrac1{u_n}\T^n\v_z\xrightarrow[n\to\infty]{}V(z)\ \text{uniformly on}\ Y
\ \text{with}\ u_n:=\T^n\v_i(0);\\ &\tag*{\faThumbsOUp}\label{faThumbsOUp}\
\sup_{n\ge 1}\|F_{n,z}\|_{\text{\tt\tiny Lip(Y)}}<\infty
\end{align*}
\endproclaim
\demo{Proof of \ref{Pointinghand}}\
\

Write $T^n(z)=\a_n(z)+i\b_n(z)\ \&\ T^n(i)=a_n+ib_n.$ By Pommerenke's theorem \ref{Pom}
$$\tfrac{\a_n(z)-a_n}{b_n}\xrightarrow[n\to\infty]{}U(z),\ \tfrac{\b_n(z)}{b_n}\xrightarrow[n\to\infty]{}V(z).$$
\

It follows that $\b_n(z)\sim b_nV(z)$ and that
\begin{align*}\tag*{\faStarO}\label{faStarO}\a_n(z)^2+\b_n(z)^2&=(a_n+U(z)+o(b_n))^2+V(z)^2b_n^2+o(b_n^2)\\ &=
(a_n+U(z))^2+V(z)^2b_n^2+2a_no(b_n)+o(b_n^2)\\ &=a_n^2+V(z)^2b_n^2+o(a_n^2+b_n^2)
\end{align*}

Thus,
\begin{align*}\tag*{\dsaeronautical}\label{dsaeronautical}
 F_{n,z}(x)&=\tfrac1{u_n}\T^n\v_z(x)=\tfrac{\b_n(z)}{b_n}\cdot
 \tfrac{a_n^2+b_n^2}{(\a_n(z)-x)^2+\b_n(z)^2}\\ &\sim V(z)\tfrac{a_n^2+b_n^2}{a_n^2+V(z)^2b_n^2}\ \ \text{uniformly on $Y$}
 \\ & =\
 \tfrac{V(z)(\D_n+1)}{\D_n+V(z)^2}\ \ \text{with}\ \
 \D_n:=(\tfrac{a_n}{b_n})^2.
\end{align*}
In case $T$ is exact,  $V\equiv 1$ entailing
$$\tfrac1{u_n}\T^n\v_z\xrightarrow[n\to\infty]{\text{\tt\tiny uniformly on} \ Y}1=V(z).$$
Otherwise by Theorem\ \ref{Pom}(iii), \ $\D_n\to \infty$ entailing
$$\tfrac1{u_n}\T^n\v_z\xrightarrow[n\to\infty]{\text{\tt\tiny uniformly on} \ Y}V(z).\ \ \CheckedBox\ \text{\ref{Pointinghand}}$$
\demo{Proof of \ref{faThumbsOUp}}
\

By \ \ \ref{Pointinghand},\ \  $\sup_{n\ge 1}\|F_{n,z}\|_{C_{\text{\tiny lim}}(X)}<\infty\ \forall\ z\in\Bbb R^{2+}$.
\

Thus, for  $\sup_{n\ge 1}\|F_{n,z}\|_{\text{\tt\tiny Lip(Y)}}<\infty$, it suffices  that each $F_{n,z}$ is absolutely continuous with
\begin{align*}\tag*{{\small\faShip}}\label{faShip}
\sup_{n\ge 1}\|F_{n,z}^{\ \prime}\|_{L^\infty(Y,m)}<\infty. \end{align*}

\pf of\  \ref{faShip}:\ \  Uniformly in $x\in Y$:

\begin{align*}
 |F_{n,z}^\prime(x)|&=\tfrac{\b_n(z)}{\pi u_n}
 |\tfrac{d}{dx}(\tfrac1{(\a_n(z)-x)^2+\b_n(z)^2})|\\ &=
 \tfrac{\b_n(z)}{\pi u_n}\tfrac{2|\a_n(z)-x|}{((\a_n(z)-x)^2+\b_n(z)^2)^2}\\ &\sim
  \tfrac{\b_n(z)}{\pi u_n}\tfrac{2|\a_n(z)|}{(\a_n(z)^2+\b_n(z)^2)^2}\ \ \text{uniformly on}\ Y
 \\ &\sim V(z)\tfrac{2(a_n^2+b_n^2)|\a_n(z)|}{((a_n^2+V(z)^2b_n^2)^2}\\ &
 \ll \tfrac1{\sqrt{a_n^2+b_n^2}}
 \xrightarrow[n\to\infty]{}0                                                                                                                                                                                                                                                                                                                                                                                                                                                                    .\ \ \CheckedBox\ \ \text{\ref{faShip}}
 \end{align*}

\subsection{Proof of the Main Theorem}
\

We show that for $f\in C_C(\Bbb R)$,
\begin{align*}\tag*{\dsmedical}\label{dsmedical}\mu_{n,x}(\tfrac{f}{\v_i})=
\tfrac1{u_n}\T^nf\xrightarrow[n\to\infty]{\text{\scriptsize u.c.s.}}\int_\Bbb Rfd\rho
\end{align*}
  with\ $u_n:=\T^n\v_i(0),\ d\rho(t):=(1+t^2)d\varpi_\Pi$ and
$\mu_n:\Bbb R\to\mathcal{M}(\Bbb R)$ defined by
 $$\mu_{n,x}(f):=\tfrac{\T^n(f\v_i)(x)}{u_n}.$$

By Lemma \ref{transferlem}\ \  \ref{Pointinghand}\ (p.\ \pageref{Pointinghand}),\,\
\ \ \ref{dsmedical}\, holds for $f=\v_z,\ z\in\Bbb R^{2+}$.

Now fix a bounded subinterval $Y$ of $\Bbb R$.  By \ \ \ref{faThumbsOUp}, \ for any $z\in\Bbb R^{2+}$,
$$\sup_{n\ge 1}\|F_{n,z}\|_{\text{\tt\tiny Lip}(Y)}<\infty$$

where $F_{n,z}(x):=\mu_{n,x}(\tfrac{\v_z}{\v_i})=\tfrac{\T^n\v_z(x)}{u_n}.$

\

Thus  $\mathcal{K}:=\{\mu_{n}|_Y:\  n\ge 1\}$ satisfies
\ref{faLeaf} (p.\ \pageref{faLeaf})\  with respect to $\text{Span}\,J$ for
$J=\{\tfrac{\v_z}{\v_i}:\ z\in\Bbb R^{2+}\}$
\ \
and by Proposition \ref{precompacthat}\ \  is precompact in
$C(Y,\mathcal{M}(\widehat{\Bbb R}))$.

\

Next, suppose that $\nu\in \mathcal{K}'$, i.e. $\exists\ n_k\to\infty\ \&\ \nu\in C(Y,\mathcal{M}(\widehat{\Bbb R}))$ are so that $\mu_{n_k}\xrightarrow[k\to\infty]{C(Y,\mathcal{M}(\widehat{\Bbb R}))}\nu$.
\

By \ref{Pointinghand},
$$\mu_{n_k,x}(\tfrac{\v_z}{\v_i})\xrightarrow[k\to\infty]{}V(z)=\varpi_\Pi(\tfrac{\v_z}{\v_i})\ \forall\ z\in\Bbb R^{2+}$$ with the conclusion that $\nu_x=\varpi_\Pi\ \forall\ x\in Y$ and
$\mu_{n_k}\xrightarrow[k\to\infty]{C(Y,\mathcal{M}(\Bbb R))}\varpi_\Pi$.

This shows \ \ \ref{dsmedical}\, \ $\forall\ f\in  C_{\text{\tt\tiny lim}}(\Bbb R)$.
In case $T$ is exact, $\Pi\equiv i,\ \varpi_\Pi=P_i\ \&\ \rho=m$. Otherwise,
$m\perp\rho$; and $\rho\ne 0$ because $T$ is not M\"obius. This completes the proof of
 the Main \hyperlink{target:Main}{Theorem}.\ \ \CheckedBox

\section{Examples of mixing $\&$ quasi-mixing}
\prp\label{PS}\ \
Let $T:\Bbb R^{2+}\CWrighttoleftarrow$ be  a non-M\"obius, normalised, parabolic inner function
with representation
$$T(z)=z+\b+\int_\Bbb R\tfrac{d\nu(t)}{t-z}\ \text{with}\ \nu\in\mathcal{M}(\Bbb R),\ \nu\perp m\ \&\ \b\in\Bbb R.$$
\f{\rm (i)}\ \ If $\b=0$, then $(\Bbb R,m,T)$ is conservative and dual mixing with
$u_n(T)=\tfrac1{\pi\sqrt{2\nu(\Bbb R)n}}$.
\

\f{\rm (ii)} \ \ If $\b\ne 0$, then $(\Bbb R,m,T)$ is singular dual quasi-mixing:
satisfying\ \ \ref{dsliterary}\ \  with respect to $d\rho_\Pi(x):=\tfrac{d\varpi_\Pi(x)}{\v_i(x)}$; with
$u_n(T)=\tfrac{b_\infty}{\pi(n\b)^2}$ where $b_\infty:=\lim_{n\to\infty}\im\,T^n(i)<\infty$ and $\varpi_\Pi$ is the representing measure of $\Pi$
as in \ \ref{faLock}.
\endproclaim
By \cite[Theorem P]{Hamilton}, $\rho$ in (ii) is singular and $\rho\nequiv 0$
 because $T$ is not M\"obius whence $\varpi_\Pi\neq\d_\infty$.

\demo{Proof of (i)}\ \ Set $T^n(i):=a_n+ib_n$, then  $b_n\sim\sqrt{2\nu(\Bbb R)n}\ \&\
a_n=o(b_n)$ as $n\to\infty$. This is established in \cite[Theorem 6.4.1]{A1} for $\nu$ compactly supported and may be extracted from \cite[\S4]{ASW} for general $\nu$.
\

Thus $(\Bbb R,m,T)$ is conservative, hence exact (see \cite[Chapter 6]{A1}); and thus  dual mixing by the Main \hyperlink{target:Main}{Theorem} with $u_n=\v_{T^n(i)}(0)\sim\tfrac1{\pi\sqrt{2\nu(\Bbb R)n}}$.\ \ \CheckedBox\ (i)
\demo{Proof of (ii)}\ \ Here $a_{n+1}-a_n\to\b$ and $b_n\uparrow b_\infty<\infty$
 as $n\to\infty]$.
See \cite[Theorem 6.4.1]{A1} for $\nu$ compactly supported and \cite{Ivrii-parabolic}
 for general $\nu$.
 \

 It follows that $\Pi\circ T=\Pi+\tfrac{\b}{b_\infty}$ and $T$ is not exact, whence $\rho_\Pi\perp m$. The  $\rho_\Pi$-quasi mixing of
 $(\Bbb R,m,T)$ is established by the Main \hyperlink{target:Main}{Theorem} with with $\v_{T^n(i)}(0)\sim \tfrac{b_\infty}{\pi(n\b)^2}=:u_n(T)$.\ \CheckedBox\ \  (ii)

\end{document}